\documentclass[11pt]{article}
\usepackage{amsmath}
\usepackage{amssymb}
\usepackage{eucal}
\usepackage{stmaryrd}
\usepackage[usenames]{color}

\begin{document}

\baselineskip 16pt

\title{A $G$-covering subgroup system of a finite group \\
 for some classes of $\sigma$-soluble groups   \thanks{Research was supported by the NNSF  of China (No. 11771409.)} }

\author{A-Ming Liu, \  W. Guo, \\
{\small School of Science, Hainan University, Haikou, 570228, P.R. China}\\ 
{\small  E-mail: amliu@hainanu.edu.cn,  wbguo@ustc.edu.cn} \\ \\
Inna N. Safonova \\
{\small Department of Applied Mathematics and Computer Science, }\\
{\small Belarusian State University, Minsk 220030, Belarus}\\
{\small E-mail: safonova@bsu.by}\\
\\
Alexander N. Skiba\\
{\small Department of Mathematics and Technologies of Programming,}\\
{\small  Francisk Skorina Gomel State University,}\\
{\small Gomel 246019, Belarus}\\
{\small E-mail: alexander.skiba49@gmail.com}}

\date{}
\maketitle

\begin{abstract} Let ${\frak F}$ be a class of group and $G$  a finite group.  
Then a set $\Sigma $ of
subgroups of $G$ is called  a \emph{$G$-covering subgroup system} for the class
${\frak F}$  if $G\in {\frak F}$ whenever  $\Sigma \subseteq {\frak
F}$.

We prove that: {\sl If a set of subgroups  $\Sigma$ of $G$ contains at 
least one supplement to each maximal subgroup of every Sylow subgroup of 
$G$, then $\Sigma$ is a $G$-covering subgroup system for the classes of 
all $\sigma$-soluble and  all $\sigma$-nilpotent groups, and for the class of 
all $\sigma$-soluble $P\sigma T$-groups.}

This result gives  positive
answers to questions 19.87 and 19.88 from the Kourovka notebook.

\end{abstract}


\let\thefootnote\thefootnoteorig

\section{Introduction}

Throughout this paper, all groups are finite and $G$ always denotes
a finite group. Moreover,  $\mathbb{P}$ is the set of all  primes and 
$\sigma$  is some partition of  
$\Bbb{P}$, that is,  $\sigma =\{\sigma_{i} \mid
 i\in I \}$, where   $\Bbb{P}=\bigcup_{i\in I} \sigma_{i}$
 and $\sigma_{i}\cap
\sigma_{j}= \emptyset  $ for all $i\ne j$.

Before continuing, recall some concepts of the papers [1, 2, 3, 4] which play 
fundamental role in the theory of $\sigma$-properties of groups.

The group    $G$ is said to be: \emph{$\sigma$-primary} if
  $G$ is a $\sigma_{i}$-group for some $i=i(G)$; 
\emph{$\sigma$-decomposable}   or 
\emph{$\sigma$-nilpotent}    if $G=G_{1}\times \dots \times G_{n}$ 
for some $\sigma$-primary groups $G_{1}, \ldots, G_{n}$;
 \emph{$\sigma$-soluble} 
 if every chief factor of $G$ is $\sigma$-primary.

A     set  ${\cal H}$ of subgroups of $G$ is a
 \emph{complete Hall $\sigma $-set} of $G$   if
 every member $\ne 1$ of  ${\cal H}$ is a Hall $\sigma _{i}$-subgroup of $G$
 for some $\sigma _{i} \in \sigma$ and ${\cal H}$ contains exactly one Hall
 $\sigma _{i}$-subgroup of $G$ for every  $i$.

 Recall  that a  subgroup $A$ of $G$ is said to be 
    \emph{$\sigma$-permutable} in $G$  if $G$ possesses
 a complete Hall $\sigma$-set  ${\cal H}$ such that $AH^{x}=H^{x}A$ for all
  $H\in {\cal H}$ and  all $x\in G$
.

We say that  $G$ is a {\sl $P\sigma T$-group} if ${\sigma}$-permutability   
is a transitive relation in $G$, that is, if $K$ is a ${\sigma}$-permutable subgroup
 of $H$ and 
 $H$ is a ${\sigma}$-permutable subgroup of $G$, then  $K$ is a
 ${\sigma}$-permutable subgroup of $G$.  

Let ${\frak F}$ be a class of group and $G$  a finite group.  
Then a set $\Sigma $ of
subgroups of $G$ is called  a \emph{$G$-covering subgroup system} \cite{Israel}
 for the class
${\frak F}$  if $G\in {\frak F}$ whenever  $\Sigma \subseteq {\frak
F}$.

In this paper, we prove the following 

{\bf Theorem A.}  {\sl Suppose that  a set of subgroups  $\Sigma$ of $G$ contains at 
least one supplement to each maximal subgroup of every Sylow subgroup of 
$G$. Then $ \Sigma  $ is a $G$-covering subgroup system for
any class  $\frak F$ in the following list:  }

(i) {\sl $\frak F$ is the class of all $\sigma$-soluble groups. }

(ii) {\sl $\frak F$ is the class of all $\sigma$-nilpotent groups. }

(iii) {\sl $\frak F$ is the class of all $\sigma$-soluble  $P\sigma 
T$-groups. }

The theory of  $P\sigma 
T$-groups was built in the works \cite{1, 2, 3}. Theorem A  gives  positive
answers to questions 19.87 and 19.88 from the Kourovka notebook \cite{30} and, also, 
  allows us to give the following new characterization of $\sigma$-soluble $P\sigma 
T$-groups.

{\bf Corollary 1.1.} {\sl $G$ is a $\sigma$-soluble $P\sigma 
T$-group if and only if each maximal subgroup of every Sylow subgroup of $G$ has a supplement 
$T$ in $G$ such that $T$ is a  $\sigma$-soluble $P\sigma 
T$-group. }

In the classical case
 when $\sigma =\sigma ^{1}=\{\{2\}, \{3\}, \ldots 
\}$: $G$  is $\sigma ^{1}$-soluble (respectively $\sigma ^{1}$-nilpotent) if
 and only if $G$ is soluble 
(respectively nilpotent); 
${\sigma} ^{1}$-permutable subgroups are   also called  
 \emph{$S$-permutable} \cite{prod}; in this case  a  $P\sigma 
T$-group is  also called a \emph{$PST$-group} \cite{prod}.

A significant place to the theory of $PST$-groups is given in the book \cite{30}.
 From Theorem A we get also the following result in this line researches.

{\bf Corollary 1.2.} {\sl $G$ is a soluble $PST$-group if and only if
 each maximal subgroup of every Sylow subgroup of $G$ has a supplement 
$T$ in $G$ such that $T$ is a  soluble $PST$-group. }

\section{Basic lemmas}

If
 $n$ is an integer, the symbol $\pi (n)$ denotes
 the set of all primes dividing $n$; as usual,  $\pi (G)=\pi (|G|)$, the set of all
  primes dividing the order of $G$.   
 $G$ is said to be  a $D_{\pi}$-group if $G$ possesses a Hall 
$\pi$-subgroup $E$ and every  $\pi$-subgroup of $G$ is contained in some 
conjugate of $E$.

By the analogy with the notation   $\pi (n)$, we write  $\sigma (n)$ to denote 
the set  $\{\sigma_{i} |\sigma_{i}\cap \pi (n)\ne 
 \emptyset  \}$;   $\sigma (G)=\sigma (|G|)$.  $G$ is said to be: a 
 \emph{$\sigma$-full group
 of Sylow type}  \cite{1} if every subgroup $E$ of $G$
 is a $D_{\sigma _{i}}$-group for every
$\sigma _{i}\in \sigma (E)$.  
 
{\bf Lemma 2.1} (See Theorem A \cite{20}). {\sl
 Every  $\sigma$-soluble group  is a $\sigma$-full group of Sylow type.}

{\bf Lemma 2.2} (Theorem 1 in \cite{000}). {\sl $G$ is $\pi$-separable if 
and only if}

(i) {\sl $G$ has a Hall $\pi$-subgroup and a Hall $\pi'$-subgroup;}

(ii) {\sl $G$ has a Hall $\pi\cup \{p\}$-subgroup and a Hall $\pi'\cup \{q\}$-subgroup 
for all $p\in \pi'$ and $q\in \pi$. }

{\bf Lemma 2.3 } (See  Corollary 2.4 and Lemma 2.5  in \cite{1}).  {\sl  
The class   of all  $\sigma$-nilpotent groups
 ${\mathfrak{N}}_{\sigma}$               is closed under taking  
products of normal subgroups, homomorphic images and  subgroups. Moreover, if  $E$ is a normal 
subgroup of $G$ and  $E/(E\cap \Phi (G))$ is $\sigma$-nilpotent, then 
$E$ is $\sigma$-nilpotent.    }

In view of Lemma 2.3, the class ${\mathfrak{N}}_{\sigma}$, of all 
$\sigma$-nilpotent groups, is a hereditary saturated formation and so
 from Proposition 2.2.8  in \cite{15} we get the following

{\bf Lemma 2.4} (See Proposition 2.2.8  in \cite{15}).   {\sl If 
$N$ is a normal subgroup of $G$, then
 $(G/N)^{{\frak{N}}_{\sigma}}=G^{{\frak{N}}_{\sigma}}N/N.$  }

In this lemma,     $G^{\frak{N_{\sigma}}}$  denotes the \emph{$\sigma$-nilpotent
 residual} of $G$,
 that is,  the intersection of all normal subgroups $N$ of $G$ with 
$\sigma$-nilpotent quotient $G/N$.

{\bf Lemma 2.5} (See Theorem A in \cite{2}).  {\sl 
 If  $G$ is a  $\sigma$-soluble 
 $P\sigma 
T$-group and  $D=G^{\frak{N_{\sigma}}}$, then  
    the following conditions hold:} 

(i) {\sl $G=D\rtimes M$, where $D$   is an abelian  Hall
 subgroup of $G$ of odd order, $M$ is $\sigma$-nilpotent  and  every element of $G$ induces a
 power automorphism in $D$;  }

(ii) {\sl  $O_{\sigma _{i}}(D)$ has 
a normal complement in a Hall $\sigma _{i}$-subgroup of $G$ for all $i$.}

{\sl Conversely, if  Conditions (i) and (ii) hold for  some subgroups $D$ and $M$ of
 $G$, then $G$ is  a $P\sigma 
T$-group.}

\section{Proof of Theorem A}

\

{\bf Proof of Theorem A. } Assume that this theorem is false.   
We can assume without loss of generality that 
$\sigma (G)=\{\sigma _{1}, \sigma _{2},  \ldots , \sigma _{t}\}$.

(I) {\sl $G$ is not $\sigma$-nilpotent. Hence $t > 1$ and
 $D:=G^{{\mathfrak{N}}_{\sigma}}\ne 1$. } 

Indeed, assume that $G$ is $\sigma$-nilpotent. Then $G$ is 
$\sigma$-soluble. Hence Statements (i) and (ii) hold for $G$. Moreover, in 
this case 
 for every $i$ the product $H_{i}$,   of all normal $\sigma _{i}$-subgroups of $G$,
 is the unique normal   Hall $\sigma _{i}$-subgroup of $G$ and $G=H_{1}\times H_{2}
 \times\cdots \times H_{t}$. Hence every subgroup of $G$ is $\sigma$-permutable in $G$.
 Thus Statement (iii) also holds for $G$, contrary to our assumption on 
$G$. Hence (I)  hods. 
  
(i) Assume that this assertion is  false and let $G$ be a counterexample 
of minimal order. 

(*)  {\sl $G$ has no  
non-identity normal   $\sigma$-primary subgroups.}

Assume that $G$ has a  minimal normal $\sigma$-primary subgroup, $R$ say. 

Let $P/R$ be any non-identity   Sylow  subgroup of $G/R$. Then for some
 prime $p$ and for a 
Sylow $p$-subgroup $G_{p}$ of $G$ we have $G_{p}R/R=P/R$, so $G_{p}$ is 
non-identity. 

Now let $V/R$ be any maximal subgroup of $P/R$, that is, $|P:V|=|P/R:V/R|=p$.
 Then  
$V=R(G_{p}\cap V)$, so
 $$p=|G_{p}R:R(G_{p}\cap V)|=
(|G_{p}||R|:|G_{p}\cap R|):(|R||(G_{p}\cap V|:|R\cap G_{p}\cap V|)=|G_{p}:G_{p}\cap V|,$$  
so $G_{p}\cap V$ is a maximal subgroup of $G_{P}$. Hence $G$ has a 
$\sigma$-soluble 
 subgroup $T$ such that $(G_{p}\cap V)T=G$.

 But then $RT/R\simeq T/(T\cap R)$ is a $\sigma$-soluble 
 subgroup of $G/R$  such that  
$$(V/R)(RT/R)=(R(G_{p}\cap V)/R)(TR/R)=G/R.$$ Therefore the hypothesis holds 
for $G/R$, so $G/R$ is $\sigma$-soluble 
 by the choice of $G$. But then $G$ is$\sigma$-soluble, a contradiction. 
 Hence we have (*).

(**)  {\sl $t =2$, that is,
 $\sigma (G)=\{\sigma _{1}, \sigma _{2}\}$.  }

Assume that $t > 2$ and let $P_{i}$ be a Sylow $p_{i}$-subgroup of $G$ for some 
$p_{1}\in  \sigma _{1}\cap \pi (G)$,  $p_{2}\in  \sigma _{2}\cap \pi (G)$  and $p_{3}\in  \sigma 
_{3}\cap \pi (G)$.  Let $V_{i}$ be a maximal subgroup of $P_{i}$.  Then, by hypothesis,
  $G$ has $\sigma$-soluble subgroups $T_{1}$, $T_{2}$ and $T_{3}$ such that 
$G=V_{i}T_{i}  $ for $i=1, 2, 3.$

Let $R$ be a minimal normal subgroup of $T_{1}$.  Then $R$ is 
$\sigma$-primary, $R$ is a  $\sigma _{k}$-group say. Since
 $|G:T_{2}|=|T_{1}T_{2}:T_{2}|=|T_{1}:T_{1}\cap 
T_{2}|$ is a $p_{2}$-number  and  $|T_{1}:T_{1}\cap 
T_{3}|$ is a $p_{3}$-number, where $p_{2}\in  \sigma _{2}$  and $p_{3}\in  \sigma 
_{3}$, we have either $R\leq T_{1}\cap 
T_{2}$ or $R\leq T_{1}\cap 
T_{3}$, $R\leq T_{1}\cap 
T_{2}$ say. Hence $R^{G}=R^{T_{1}T_{2}}=R^{T_{2}}\leq T_{2}$, so $G$ has a 
non-identity normal   $\sigma$-primary subgroup, contrary to Claim (*). 
Thus (**) holds.

{\sl The final contradiction for (i). }  Let $\pi =\sigma _{1}\cap \pi (G)$.  Since 
$T_{i}$ is $\sigma$-soluble, $T_{i}$ has a Hall $\sigma _{k}$-subgroup 
for all $k$ by Lemma 2.1. Then  a Hall $\sigma _{2}$-subgroup of $T_{1}$ is a 
Hall $\pi'$subgroup of $G$ and a Hall $\sigma _{1}$-subgroup of $T_{2}$ is a 
Hall $\pi-$subgroup of $G$. 

Now we show that $G$ has a Hall $\pi\cup \{p\}$-subgroup for every 
$p\in \sigma _{2}\cap \pi (G)$. If $|\sigma _{2}\cap \pi (G)= 1$ it is 
evident. Now assume that $|\sigma _{2}\cap \pi (G)| > 1$ and
 let $q\in (\sigma _{2}\cap \pi (G))\setminus \{p\}.$ Let $V$ be a 
maximal subgroup of a Sylow $q$-subgroup $Q$ of $G$. And let $T$ be a 
$\sigma$-soluble supplement to $V$ in $G$. Then $T$ is $\pi$-separable by 
Claim (**). Hence $T$ has a Hall $\pi \cup \{p\}$-subgroup  $H$ by Lemma 
2.2.  But $|G:T|$ is a $\{q\}$-number, where $p\ne q\not \in \sigma _{1}$,
 so $H$ is a 
Hall $\pi\cup \{p\}$-subgroup of $G$.

  Similarly it can be proved that 
$G$ has a Hall $\pi'\cup \{p\}$-subgroup for all  $p\in \pi$. Therefore 
$G$ is $\pi$-separable by Lemma 2.2 and  so $G$ is $\sigma$-soluble by Claim 
(**), contrary to the choice of $G$.  Hence Statement (i) holds.

(iii) Assume that this assertion is  false and let $G$ be a counterexample 
of minimal order. Then $G$  is $\sigma$-soluble by Part (i).

(1) {\sl If $R$ is a non-identity normal subgroup of $G$, then the 
 hypothesis holds for $G/R$. Hence $G/R$ is a $\sigma$-soluble $P\sigma T$-group}
 (See the proof of Claim (*)).

(2) {\sl If $R$ is an abelian minimal normal subgroup of $G$, then $R$ is not a
 Sylow subgroup of $G$.}

Indeed, assume that $R$ is Sylow subgroup of $G$  and let $V$ be an
y maximal subgroup of $R$. Then for every supplement $T$ to $V$ in $G$ we have that
 $T\cap R$ is normal in $G$, the minimality of $R$ implies that $T=G$. 
Hence $G$ is  a $\sigma$-soluble $P\sigma T$-group, a contradiction. Hence 
(2) holds.

(3) {\sl $D$ is $\sigma$-nilpotent.}

Assume that this is falls.  Then $D$ is not $\sigma$-primary.
 Let $R$ be a minimal normal subgroup of $G$,
 so $R$ is a $\sigma _{i}$-group for some 
$i$ since $G$ is $\sigma$-soluble. 
 Moreover, from Lemmas 2.3 and 2.4 we get that 
$$(G/R)^{{\mathfrak{N}}_{\sigma}}=G^{{\mathfrak{N}}_{\sigma}}R/R=
DR/R\simeq D/(D\cap R)$$ is  a Hall $\sigma$-nilpotent subgroup of $G/R$ 
 by Claim (1).  Hence $R$ is the unique minimal normal subgroup of $G$, $R < D$
 and 
$R\nleq \Phi (G)$ since $D$ is not $\sigma$-nilpotent.
 Therefore $ C_{G}(R)\leq R$ and $D/R$ is
 a Hall subgroup of $G/R$.  Moreover, $D/R$ is not a $\sigma _{i}$-group 
since $D$ is not $\sigma$-primary. Let $p$ be a prime dividing $|D/R|$ such 
that $p\not \in  \sigma _{i}$. And let $P$ be a Sylow $p$-subgroup of $D$. 
Then $P\cap  R=1$ and $P$ is a Sylow $p$-subgroup of $G$ since $D/R$ is a 
Hall subgroup of $G/R$. 

Let  $V$ be a maximal subgroup of $P$ and $T$ a supplement to $V$ in $G$ 
such that $T$ is a $P\sigma T$-group. Then 
$T^{{\mathfrak{N}}_{\sigma}}\leq D$ and $T^{{\mathfrak{N}}_{\sigma}}$ is a Hall abelian
 subgroup
 of $T$ such that every subgroup of $T^{{\mathfrak{N}}_{\sigma}}$ is normal in $T$ by
 Lemma 2.5.
Moreover, $R\leq T$ since $|G:T|$ is a $\{p\}$-number. Hence 
$T^{{\mathfrak{N}}_{\sigma}}\cap R$ is a normal abelian Hall subgroup of 
$R$. Hence either $T^{{\mathfrak{N}}_{\sigma}}\cap R=1$  or 
$T^{{\mathfrak{N}}_{\sigma}}\cap R=R$ and so $R\leq 
T^{{\mathfrak{N}}_{\sigma}}$.  

First assume that  $T^{{\mathfrak{N}}_{\sigma}}\cap R=1$. Then 
$T^{{\mathfrak{N}}_{\sigma}}\leq C_{G}(R)$, so 
$T^{{\mathfrak{N}}_{\sigma}}=1$ and hence $T$ is $\sigma$-nilpotent. 
From $P=P\cap VT=V(P\cap T)$ it follow that $T$ is not a $\sigma 
_{i}$-group. Hence for a Hall $\sigma _{i}'$-subgroup  $E$ of $T$ we have 
$E\ne 1$ and $E\leq C_{G}(R)\leq R$, a contradiction.  Therefore $R\leq 
T^{{\mathfrak{N}}_{\sigma}}$, so $R=T^{{\mathfrak{N}}_{\sigma}}$ is a $q$-group for
 some prime 
$q\ne p$ since $C_{G}(R)\leq R$ 
and $T^{{\mathfrak{N}}_{\sigma}}$ is abelian. Let $Q$ be a Sylow 
$q$-subgroup of $T$.  Then $R=Q$ since  $R=T^{{\mathfrak{N}}_{\sigma}}$ is 
a Hall subgroup of $T$.    
Moreover, $R$ is a Sylow $q$-subgroup of $G$ since $p\ne q$ and  
$|G:T|$ is a $\{p\}$-number, contrary to Claim (2). Hence we 
have (3).

(4) {\sl $D$ is nilpotent.}

Assume that this false and let $R$ be a minimal normal subgroup of $G$.
Then $R\leq D$ and $C_{G}(R)\leq R$ and $D/R$ is a Hall subgroup of $G/R$ (see the 
proof of Claim (3)).
 Hence $D\leq O_{\sigma _{i}}(G)$ for 
some $i$ by Claim (3). Let $P$ be a Sylow $p$-subgroup of $G$,
 where $p\in \pi (G)\setminus \sigma _{i}$. Let $V$ be a maximal subgroup 
of $P$ and $T$ a supplement to $V$ in $G$ such that $T$ is a $P\sigma 
T$-group. Then $R\leq D\leq T$, so  $T^{{\mathfrak{N}}_{\sigma}}\cap R$
 is a normal abelian Hall subgroup of 
$R$. Hence  $R\leq 
T^{{\mathfrak{N}}_{\sigma}}$ (see the proof of Claim (*)). On the other 
hand, $T^{{\mathfrak{N}}_{\sigma}}\leq D$. Therefore 
$R=T^{{\mathfrak{N}}_{\sigma}}$ is a Sylow $q$-subgroup of $G$ for some 
$q\ne p$, contrary to Claim (2).

(5) {\sl  $D$ is a Hall subgroup of $G$. }

 Suppose
that this is false and let $P$ be a  Sylow $p$-subgroup of $D$ such
that $1 < P < G_{p}$, where $G_{p}\in \text{Syl}_{p}(G)$.  We can assume 
without loss of generality that $G_{p}\leq H_{1}$.

(a)  {\sl    $D=P$ is  a minimal normal subgroup of $G$. }

Let $R$ be a minimal normal subgroup of $G$ contained in $D$. 
 Since
 $D$ is  nilpotent by Claim (4),   $R$ is a $q$-group    for some prime   
$q$. Moreover, 
$D/R=(G/R)^{\mathfrak{N}_{\sigma}}$  is a Hall subgroup of $G/R$ by
Claim (1) and Lemma 2.3.  Suppose that  $PR/R \ne 1$. Then  $PR/R \in \text{Syl}_{p}(G/R)$. 
If $q\ne p$, then    $P \in \text{Syl}_{p}(G)$. This contradicts the fact 
that $P < G_{p}$.  Hence $q=p$, so $R\leq P$ and therefore $P/R \in 
\text{Syl}_{p}(G/R)$ and we again get  that 
$P \in \text{Syl}_{p}(G)$. This contradiction shows that  $PR/R=1$, which implies that 
  $R=P$  is the unique minimal normal subgroup of $G$ contained in $D$.
 Since $D$ is nilpotent by Claim (4),
 a $p'$-complement $E$ of $D$ is characteristic in 
$D$ and so it is normal in $G$. Hence $E=1$, which implies that $R=D=P$.

(b) {\sl $D\nleq \Phi (G)$.    Hence for some maximal subgroup
 $M$ of $G$ we have $G=D\rtimes M$.  }

(c) {\sl If $G$ has a minimal normal subgroup $L\ne D$, then $G_{p}=D\times (L\cap G_{p})$.
  Hence $O_{p'}(G)=1$. }

Indeed, $DL/L\simeq D$ is a Hall 
subgroup of $G/L$ by Claim (1). Hence  $G_{p}L/L=RL/L$, so $G_{p}=D\times (L\cap G_{p})$.
 Thus  $O_{p'}(G)=1$ since $D < G_{p}$ by Claim (a).

(d)  {\sl   $V=C_{G}(D)\cap M$ is a  normal subgroup of $G$ and 
 $C_{G}(D)=D\times V \leq H_{1}$.  }

In view of  Claim  (b),  $C_{G}(D)=D\times V$, where $V=C_{G}(D)\cap M$ 
is a normal  subgroup of $G$. By Claim (a), $V\cap D=1$ and 
hence $V\simeq DV/D$ is $\sigma $-nilpotent by Lemma 2.2.  Let $W$ be a $\sigma 
_{1}$-complement of $V$. Then $W$  is characteristic in $V$ and so it is normal 
in $G$.    Therefore we have  (d) by Claim (c).

{\sl The final contradiction for (5).} Let $Q$ be a Sylow $q$-subgroup of $G$, where 
$q\in \pi (G)\setminus \pi (H_{1})$. Let  $V$ be a maximal subgroup of $P$ and $T$ a supplement to $V$ in $G$ 
such that $T$ is a $P\sigma T$-group. Then 
$T^{{\mathfrak{N}}_{\sigma}}\leq D$ and $T^{{\mathfrak{N}}_{\sigma}}$ is a Hall abelian
 subgroup   of $T$. Then $D$ is not a Sylow $q$-subgroup of $T$ and so 
$T^{{\mathfrak{N}}_{\sigma}}=1$, which implies that $T$ is $\sigma 
$-nilpotent. But  then for a Sylow $q$-subgroup $T_{q}$ 
of $T$ we have   $1 < T_{q} \leq C_{G}(D)\leq H_{1}$, a contradiction.

(6) {\sl Every subgroup $H$ of $D$ is normal in $G$. Hence every element of
 $G$ induces a power automorphism in $D$. }

Since $D$ is  nilpotent by Claim (4), it is enough
 to consider   the case  when $H\leq O_{p}(D)$ for some $p \pi (D)$.

Let $R$ be any Sylow $r$-subgroup of $G$, where $r\not \in \pi (D)$. 
 Let $V_{1}, V_{2}, \ldots , V_{t}$ be the set 
of all maximal subgroups of $R$. Let $T_{i}$ be a supplement to $V_{i}$ in $G$ such that 
$T_{i}$ is a $P\sigma T$-group with $D_{i}=T^{\mathfrak{N}}$.

Since $G=V_{i}T_{i}$, $R=V_{i}(T_{i}\cap R)$. Hence for some $a_{i}\in T_{i}\cap R$
 we have 
 we have $a_{i}\not\in V_{i}$. We show that $a_{i}\in N_{G}(H)$.

First observe that $ D\leq T_{i}$ since $|G:T|$ is a $q$-number,
 where $r\not \in \pi (D)$. 
 Moreover,  $D_{i}\leq D$. But $D_{i}$ is a Hall
 subgroup of $T_{i}$ and every subgroup of $D_{i}$ is normal in $T_{i}$,
 so $D_{i}$ is a Hall subgroup of $D$. So either $H\leq O_{p}(D)\leq 
D_{i}$ or  $O_{p}(D)\cap D_{i}=1$.   In the former case we have $a_{i}\in 
N_{G}(H)$ since every subgroup of $D_{i}$ is normal in $T_{I}$. Now assume 
that $O_{p}(D)\cap D_{i}=1$, so $D_{i}\cap O_{p}(D)\langle a_{i}\rangle=1$ 
since $D_{i}\leq D$  and $r\not \in \pi (D)$,  so $ O_{p}(D)\langle a_{i}\rangle\simeq D_{i}O_{p}(D)\langle 
a_{i}\rangle/D_{i}$ is $\sigma$-nilpotent. Hence  $[O_{p}(D), a_{i}]=1$, 
so $a_{i}\in N_{G}(H)$.

  Let $V=\langle  a_{1}, a_{2}, \ldots , 
a_{i}\rangle $. Then $V\leq N_{G}(H)$.  
 Moreover, if $V < R$, then 
 for some  ${i}$ we have $V\leq V_{i}$.
 But then  $a_{i}\not\in V_{i}$ and 
$a_{i}\in V\leq V_{i}\leq V_{i}$, a contradiction. Therefore $V=R\leq 
N_{G}(H)$. Hence $R^{G} \leq N_{G}(H)$. 
Therefore  $E^{G}\leq N_{G}(H)$, where $E$ is a Hall $\pi (D)'$-subgroup 
of $G$. But then $E^{G}D/E^{G}\simeq D/(D\cap E^{G})$ is nilpotent, so 
$D\leq E^{G}$ and hence $G=GE=E^{G}$. Hence we have (6).

 (7) {\sl  If  $p$ is a  prime such that $(p-1, |G|)=1$, then  $p$
does not divide $|D|$. In particular,  $|D|$ is odd. }

Assume that this is false.
 Then, by Claim (4),  $D$ has a maximal subgroup $E$ such that
$|D:E|=p$ and  $E$ is normal in $G$. It follows that  $C_{G}(D/E)=G$ since $(p-1, 
|G|)=1$.   
Since
$D$ is a Hall subgroup of $G$, it has a complement $M$ in $G$. 
Hence 
$G/E=(D/E)\times (ME/E)$, where
$ME/E\simeq M\simeq G/D$ is $\sigma$-nilpotent. Therefore $G/E$ is
$\sigma$-nilpotent. But then $D\leq E$, a contradiction. Hence $p$
does not divide $|D|$. In particular, $|D|$  is odd.

(8)  {\sl  $D$ is abelian.}

In view of Claim 
(5), $D$ is a Dedekind group.  Hence $D$ is abelian since $|D|$ is  odd  by Claim (7).  

 From Claims (5)--(8)
 we get that $G$ is $\sigma$-soluble $P\sigma T$-group, contrary to the 
choice of $G$. Hence Statement  (iii) holds.

(ii) Assume that this assertion is  false and let $G$ be a counterexample 
of minimal order. Then $G$ is a $\sigma$-soluble $P\sigma T$-group by Part 
(iii) since every $\sigma$-nilpotent group is a $\sigma$-soluble $P\sigma T$-group. 
Then $G^{{\mathfrak{N}}_{\sigma}}$ is a Hall subgroup of $G$ of odd order and every 
subgroup of $G^{{\mathfrak{N}}_{\sigma}}$ is normal in $G$ by lemma 2.5. Moreover, the 
hypothesis holds on  $G/R$ for every minimal normal subgroup $R$ of $G$ 
and hence 
 $G/R$ is $\sigma$-nilpotent by the choice of $G$, so $R=G^{{\mathfrak{N}}_{\sigma}}$
 is a group of prime order $p$ for some prime $p$ and $R$ is a Sylow $p$-subgroup of $G$.
But then the maximal $V$ subgroup of $R$ is identity and so $G$ is
 the unique supplement to  $V$ in $G$, so $H$ is $\sigma$-nilpotent, a contradiction. Therefore
Statement (ii) holds.    

The theorem is proved.

\end{document}